\documentclass[11pt]{amsart}
\usepackage{latexsym}
\usepackage{amssymb,amsmath}
\theoremstyle{plain}
\newtheorem{thrm}{Theorem}[section]

\def\doublespaced{\baselineskip=\normalbaselineskip
    \multiply\baselineskip by 2}
\def\doublespace{\doublespaced}
\def\singlespaced{\baselineskip=\normalbaselineskip}
\def\singlespace{\singlespaced}

\doublespace

\newcommand{\rn}[1]{\mathbb{R}^{#1}}

\newcommand{\beq}{\begin{equation}}
\newcommand{\bea}[1]{\begin{array}{#1} }
\newcommand{\eeq}{ \end{equation}}
\newcommand{\ea}{ \end{array}}
\newcommand{\ep}{\epsilon}
\newcommand{\es}{\emptyset}

\newcommand{\al}{\alpha}

\newcommand{\de}{\delta}
\newcommand{\ds}{\displaystyle}
\newcommand{\ts}{\textstyle}
\newcommand{\rar}{\mbox{$\rightarrow$}}

\newcommand{\ran}{\rangle}
\newcommand{\lan}{\langle}

\newcommand{\la}{\lambda}

\newcommand{\ar}{\partial}
\newcommand{\si}{\sigma}

\newcommand{\Om}{\Omega}

\newcommand{\ph}{\phi}
\newcommand{\he}{\theta}

\newcommand{\hs}[1]{\mbox{$ \hspace{#1}$}}

\newcommand{\sem}{\setminus}
\newcommand{\ze}{\zeta}

\newcommand{\De}{\Delta}
\newcommand{\ti}{\tilde}

\def\singlespace{\singlespaced}
\singlespace
\begin{document}
\title[Gradient bounds  etc]{Gradient bounds for p-harmonic systems with vanishing  Neumann data in a convex domain}
\author{Agnid Banerjee}
 \email[Agnid Banerjee]{banerja@math.purdue.edu}\thanks{Banerjee was supported in part  by N. Garofalo's Purdue Research Foundation Grant, 2012 and in part  by N. Garofalo's  NSF Grant DMS-1001317}
\address{Department of Mathematics\\Purdue University \\
West Lafayette, IN 47907}
     \author{John L. Lewis}
\email{johnl@uky.edu} \thanks{
Lewis was partially supported
by DMS-0900291}
\address{Department of Mathematics, University of Kentucky\linebreak \\  Lexington, KY  40506-0027, USA}
\begin{abstract}
\noindent Let $ \ti \Om $ be a  bounded convex  domain in Euclidean $ n $ space, $ \hat x \in \ar \ti \Om, $ and $ r  > 0. $    Let   $ \ti u = ( \ti u^1, \ti u^2, \dots, \ti u^N ) $  be a weak solution  to
\[   \nabla \cdot \left ( |\nabla \ti u |^{p-2}  \nabla \ti u  \right) = 0 \mbox{ in }
 \ti \Om \cap B ( \hat x, 4 r ) \mbox{ with }  |\nabla \ti u|^{p-2} \,  \ti u_\nu  = 0  \mbox{ on } \ar \ti \Om \cap B  ( \hat x, 4 r ).  \]  We show  that  sub solution type arguments for certain  uniformly elliptic systems can be used to  deduce that
   $ | \nabla \ti u | $  is  bounded in $ \ti \Om \cap B ( \hat x, r )$
   with constants depending only on $ n, p, N. $
 and $ \frac{r^n}{ | \ti \Om \cap B ( \hat x, r ) |}.$      Our argument replaces an argument based on level sets in recent
     important work of  [CM], [CM1], [GS], [M], [M1], involving similar problems.
       \\


\end{abstract}
\maketitle

\section{Introduction}
Let  $ x = ( x_1, \dots, x_n ) $ denote points in Euclidean $n$ space, $ \rn{n}, $  and  let  $ \lan \cdot, \cdot \ran $  denote the standard inner product on
$ \rn{n}. $  let  $ | x | = \lan x, x \ran^{1/2} $ denote the norm of $ x $ and set  $ B ( x, \rho  ) = \{ y : | y - x |< \rho  \} $ when $ \rho > 0. $  Given an open set $ O \subset  \rn{n}, $  let $ C_0^\infty ( O ) $ denote infinitely differentiable functions with compact support in $ O. $  If  $ 1 \leq   q  \leq  \infty, $  let  $ W^{1, q } ( O ) $ denote the Sobolev space of functions $ g $  with distributional  derivatives $ g_{x_i}, 1 \leq i \leq n, $  and norm,
\[ \| g \|_{W^{1,q} ( O )}  =  \| g \|_{L^q (O)} +   \| |\nabla g| \|_{L^q (O)}  <  \infty.  \]
Here  $ \nabla  g $ denotes the gradient of $ g $ and  $ \| \cdot \|_{L^q(O)}  $ is  the usual  Lebesque $ q $ norm relative to $ O. $                                                 If $ E, F \subset  \rn{n} $ let  $ \hat H ( E, F ) $ denote the  Hausdorff distance between the sets
$ E,  F,  $  and let   $|E|$ denote the Lebesgue $n$ measure of  $ E $ whenever $ E $ is measurable.
  Throughout this paper we assume that $ \ti \Om \subset \rn{n} $ is a  bounded convex domain  and for  given $ \hat x \in  \ti \Om, 1 < p < \infty,  r > 0, $ that $ \ti u = (\ti u^1, \dots, \ti u^N ) :  \ti \Om \rar
\rn{N} $ is a weak solution to  the $p$ Laplace systems equation in
$ \ti \Om \cap B ( \hat x,4 r ) $ with vanishing Neumann data on
$ \ar \ti \Om \cap B ( \hat x, 4 r ). $  That is,
$ \ti u^k  \in  W^{1,p} (\ti \Om \cap B ( \hat x, 4r ) ), 1 \leq k \leq N,  $
 and if $ \ph = ( \ph^1, \dots, \ph^N), $
\beq  \tag{1.0}  \int_{\ti \Om \cap B ( \hat x, 4r ) }  | \nabla \ti u |^{p-2}
\, \lan\nabla \ti u,  \nabla \ph \ran  dx = 0 \eeq
whenever   $ \ph^i \in  W^{1,p} (\ti \Om \cap B ( \hat x, 4 r ) ) $  and $ \ph^i, 1 \leq i \leq N, $   vanishes outside a set  whose closure is a  compact subset of
$ B ( \hat x,4 r ). $  In the above display,   $ \nabla \ti u$ is the $ nN $  tuple,    $ ( \nabla \ti u^1, \dots, \nabla \ti u^N ) $  while the inner product is relative to  $  \rn{nN}. $  In this note we show that
\begin{thrm}  Let $ \ti u, \hat x, r, \ti \Om, p $ be as above. There exists
$ C$ depending only on $ n, p,  N,   r^n/|\Om \cap B ( \hat x, r )|, $ such that for every $ y  \in \ti \Om \cap B ( \hat x, r ), $
   \[  | \nabla \ti u |^p (y)
 \leq  C   \,   r^{-n} \,
  \int_{ \ti \Om \cap B ( \hat x, 4 r ) } | \nabla \ti u |^p  \, dx  \]
\end{thrm}
We note that  in [CM], [CM1], the authors studied  weak solutions to  quasilinear elliptic equations and systems of  the form
$ \nabla \cdot (  a ( |\nabla \ti u| ) \nabla \ti u  )  =  f (x) $  in a  convex domain $\ti \Om $ under  both   Dirichlet  and  Neumann boundary conditions:
$ \ti u  \equiv 0,   \ti u_\nu  \equiv 0 $ respectively on  $ \ar \ti \Om. $ It is not obvious to us that their nearly endpoint global type results   imply   Theorem 1.1.   Our proof though  uses  some of  the same  arguments as in these papers, including  fundamental use of  an  inequality in [G]    for convex domains.  Also as in these papers, we first prove  Theorem 1.1 for a weak solution to a  related PDE  in  a smooth convex domain   and  then use a  limiting  argument to get Theorem 1.1.  However the  above  authors  use the   level sets of $ |\nabla u | $ and perform some rather involved calculations on these sets in order to obtain  their results.  In contrast  we  use well known sub solution type arguments  for  $ | \nabla \ti u | $  (i.e, Moser iteration)  to get Theorem  1.1.

\medskip
\noindent
\textbf{Acknowledgment:}We would like to thank Prof. N Garofalo for bringing the reference [CM1] to our attention and showing interest in our work.

\section{Proof of  Theorem 1.1 } To begin the proof
of  Theorem 1.1  let $ c $ be a positive constant, not necessarily the same at each occurrence,  which unless otherwise stated may only depend on $ n, p, N.  $  We note  that the $ p $  Laplace equation is invariant under dilations,  translations, and scaling.  Thus it suffices to prove
Theorem 1.1
  when
  \beq \tag{2.0}  r = 1, \hat x = 0,   \mbox{ and }   \int_{ \ti \Om
  \cap B ( 0, 4 )} \,    | \nabla \ti u |^p \, dx  = 1, \eeq   as we can then transfer back  to the general case  using
    the transformation  $ y =
     \hat x   +  r  x, $                                                                    after multiplying $ \ti u $  by an appropriate constant.  We continue under assumption  (2.0).  We may also assume that
     $ \ti  \Om  \cap \ar  B (0, t )  \not = \es $ whenever $ t < 4 $  since otherwise it follows from (1.0) with  $ \ph =$ a suitable extension of $ \ti u $ to  $ B ( 0, 4 ) $ that  $ \ti u \equiv $  constant.   Let    $ \rho,    3  < \rho  <  4 , $ be such that
 \beq  \tag{2.1}   \int_ {\ti \Om  \cap \ar B ( 0,  \rho )}  | \nabla \ti u |^p  \, d \si   \leq   c    \eeq  where $ \si $ denotes  Hausdorff $ n - 1 $ measure = surface area.  Existence of
   $ \rho, c $  follows from (2.0), writing the integral in polar  coordinates. and  the usual weak type estimates. Moreover, (2.1) holds for $\rho$ belonging to a set of positive measure in $[3,4]$.  Let $\eta(x)=\gamma(|x|)$ where $\gamma \equiv -1$  in $[0, \rho- \delta]$, linear in $[\rho-\delta, \rho]$ and   $ \gamma (t) = 0$ for all $t\geq \rho$. Using this $\eta$ as a test function in (1.0) and letting $\delta \to 0$, we deduce from the Lebesgue Differentiation theorem, that we  may assume
      \beq  \tag{2.2}   \int_ { \ti \Om  \cap  \ar B ( 0,  \rho )}
         |\nabla \ti u |^{p-2} \ti u_\nu  \, d \si  = 0   \eeq
   where $ \nu (x) $ denotes the outer    unit normal to
   $ \ti \Om \cap B ( 0, \rho ) $ at $ x \in \ti \Om \cap \ar B ( 0,  \rho) $ and $ \ti u_\nu  = ( \ti u^1_\nu, \dots, \ti u_\nu^N ) $ with
   $\ti u^k_\nu (x) =  \lan \nabla \ti u^k (x), \nu (x) \ran, 1 \leq k \leq N . $
      Next given  $ \ep > 0 ,  0 < \ep < 1/2, $  let $ \Om = \Om ( \ep ) $ be a convex domain with
   \beq \tag{2.3}    \ti \Om \subset \Om, \, \,
    \, \,  \hat H ( \ar \Om, \ar \ti \Om ) <  \ep, \, \, \mbox{ and }  \, \, \ar \Om  \in C^{\infty} .  \eeq
   Let  $ \rho  $ be as in  (2.2)   and let $  f = (f^1, \dots, f^N) $ be defined by  $ f^k (x)  = 0, 1 \leq k \leq N, $   when
    $ x \in  ( \Om \sem \ti \Om )  \cap \ar B ( 0,  \rho )  $
    while $ f  =  | \nabla \ti u|^{p-2}  \ti u_\nu,  $ otherwise on  $ \Om \cap
     \ar B ( 0, \rho). $  Note from (2.2) that $ f $ satisfies the compatibility condition :
\beq  \tag{2.4}   \int_{\Om \cap \ar B ( 0,  \rho )} f^k  d \si  = 0,  1 \leq k \leq N. \eeq
Using   (2.4) along with well known trace theorems  and variational methods, we deduce the existence of a unique   $ u = u ( \cdot, \ep ) : \Om \cap B ( 0, \rho) \rar R^N $   with \\ $ u^k  \in    W^{1,p} ( \Om \cap B ( 0, \rho ) ), 1 \leq k \leq N,  $  satisfying  $ \int_{\Om \cap B ( 0, \rho)} \, u \, dx = 0 $ and
\beq \tag{2.5}  \int_{ \Om \cap B ( 0, \rho ) }  ( \ep + |\nabla u |^2 )^{p/2-1}   \lan \nabla u, \nabla v \ran \, dx \, = \, \int_{\Om \cap \ar B ( 0, \rho ) }  \, \lan f,  v \ran \,  dx \eeq
whenever $ v = (v^1, \dots, v^N) $ with $ v^k \in W^{1,p} ( \Om \cap  B ( 0, \rho) ) . $  In (2.5)  we have  used $ v $ on $ \Om \cap \ar B ( 0,  \rho ) $ to denote the trace of $ v.$    More specifically,  $ u $ is the minimum of the
 functional

 \[ F_{\ep}={\ds \frac{1}{p}\int_{\Om \cap B(0,\rho)} (\ep + |Dw|^2)^{p/2}dx} -{\ds \int_{\partial B(0, \rho) \cap \Om} <f, w>} \]
 among all $w \in W^{1,p}$ such that $ \int_{\Om \cap B(0, \rho} w=0. $  Existence follows from  lower semicontinuity of $F_{\ep}$, (2.4),   and compactness of the trace operator.  Choosing $ v  = u $ in (2.5) and using (2.4), Poincar\'{e}'s inequality,  (2.1), we deduce that
 \beq \tag{2.6}  \int_{ B ( 0,  \rho ) }  ( \ep + |\nabla u |^2 )^{p/2-1} |\nabla u |^2 \, dx  \leq  c .  \eeq
 Next we claim that
   \beq  \tag{2.7}   u^k    \in C^{\infty}(\overline{ \Om  \cap B ( 0, 2 )} ) \mbox{ for }  1 \leq k \leq N \mbox{ and } 1 \leq m \leq n. \eeq
    Constants in (2.7) however may depend on $ \ep, p, n, N $ and the smoothness of  $ \ar \Om $ in (2.3). We sketch the proof of (2.7) in the appendix to this paper for the reader's convenience  since we have not been able to find a suitable reference.  Our proof uses a   reflection type argument, after  straightening  $ \ar \Om  $ in an appropriate way, to first  get $C^{1,\alpha}$ regularity  in the closure of $\Om \cap B(0,2).$  After that  we  use  Schauder type arguments  to bootstrap.
     as in [ADN] , [ADN1].

       From (2.5) and (2.7) we see that  $ u $ is a strong solution to
  \beq \tag{2.8}  \nabla \cdot \left( (\ep + | \nabla u |^2 )^{p/2 - 1} \nabla u  \right)  = 0  \eeq
  in $ \Om \cap B ( 0, \rho ), $  where $ \nabla \cdot $ denotes the divergence operator.
  Let   $  \la  \in \rn{n} $ with $ | \la | = 1. $  Let $ u_\la =  ( \nabla u^1_\la, \dots, \nabla u^N_\la ) $ denote the directional derivative of $ u $ in the direction $ \la. $
Differentiating (2.8) with respect to $ \la $ we get for fixed  $ l,  1 \leq l \leq N, $  that  $ u^m_\la $ is a solution in  $ \Om \cap B(0, \rho ) $ to  \beq \tag{2.9}     \sum^N_{
m=1}  \sum_{i,j=1}^n  \frac{\ar}{\ar x_i}  \left(  b^{l m}_{i j} {\ts \frac{\ar}{\ar x_j}} (  u^m _{\la} ) \right) = 0  \eeq
where     \beq \tag{2.10}  b^{l m}_{i j} (x)   = (  \ep + |\nabla u |^2  )^{p/2 -2}  \left(  (p-2)  u_{x_j}^m \, u_{x_i}^l   +  \de_{ij}^{l m } ( \ep  +    | \nabla u |^2  )   \right) (x)  \eeq
 when  $ x \in \Om \cap B (0, \rho). $
Here $ \de^{l m }_{i j}  = 1 $ if  $ i = j, l = m, $  and is  zero otherwise. We note that   if   $ \xi =  ( \xi_{k}^{\nu} ) $ is an  $ n \times N $ matrix,   $ | \xi |^2 = \lan \xi, \xi \ran $ where the inner product is relative to $ \rn{nN} $ and  $ h = ( \ep + |\nabla u |^2 )^{1/2}, $  then
at  $ x \in \Om \cap B ( 0, \rho) $ \beq \tag{2.11}  \min ( p - 1, 1 ) | \xi |^2  h^{p-2}   \leq   \sum_{m,l = 1}^N
\sum_{i, j = 1}^n   b_{ij}^{lm}  \xi_{i}^{l}   \xi_{j}^{m}   \leq   \max ( p - 1, 1 ) |\xi |^2 h^{p-2}  .  \eeq
For  $ 1 \leq i \leq j \leq n  $ let  \beq \tag{2.12}   c_{ij} (x)  = ( \ep + | \nabla u|^2 (x) )^{- 1} [ (p-2)  ( \sum_{\nu = 1}^N  u^l_{x_i} \, u^l_{x_j} ) + \de_{ij} ( \ep + | \nabla u |^2 ) \, \, ] \eeq
where  $ \de_{ij} $ is the  Kronecker $ \de $  and let  $ L  $ be the differential operator,
\[    \sum_{i,j = 1}^n \frac{\ar}{\ar x_i} \left( c_{ij} (x)  \frac{\ar}{\ar x_j} \, \right) \, . \]
          We observe that if  $ \mu = ( \mu_1, \dots, \mu_n ) \in \rn{n}, $ then
 \beq \tag{2.13}  \min (  p - 1, 1 )  | \mu |^2   \, \leq \,
 \sum_{i,j=1}^n  c_{ij} ( x ) \, \mu_i  \mu_j  \,  \leq \max ( p-1 , 1)
 \, | \mu |^2  \, . \eeq

   We use (2.12), (2.13)  to  give a  proof of  a well known sub solution inequality  for  $ L ( h^q ). $   To this end let $ e_k $ be the  point in $ \rn{n}$  with a  one in the $k$ th position and zeroes elsewhere.  If  $ q \geq p $  at  $ x \in \Om \cap B (0, \rho ), $ we calculate
      \beq  \bea{ll}  \tag{2.14}
   L ( h^q )  & = {\ds  \sum_{i,j = 1}^n  ( c_{ij} ( h^q )_{x_j}\, )_{x_i}}    \\ \\
    & = q (q-p)   h^{q-2}   { \ds \sum_{i,j=1}^n  \, c_{ij} } \, h_{x_i} h_{x_j}  +  (q/p) h^{q-p}  L (h^p).
    \ea  \eeq
    Moreover,  \beq   \tag{2.15}
   L ( h^p )   =
      {\ds p \sum_{ i, j, k =1}^n    \sum_{m=1}^N }  (  h^{p-4} u^l_{x_j} u^l_{x_i}  u^m_{x_k}   u^m_{x_k x_j}  )_{x_i} +   \De ( h^p)  =  S_1 + \De ( h^p) \,      \eeq  where $ \De $ is the Laplacian.    Thus
  \beq  \bea{l}  \tag{2.16}  S_1 =  p(p-2)    {\ds  \sum_{ i, j, k =1}^n   \sum_{m, l =1}^N }
    h^{p-4}  u^l_{x_j x_i}   u_{x_i}^l  u^m_{x_k}   u^m_{x_k x_j}   \\ \\
    \hs{.7in}  +   p(p-2)    {\ds  \sum_{ i, j, k =1}^n   \sum_{m, l =1}^N }       u^l_{x_j } (  h^{p-4}  u_{x_i}^l  u^m_{x_k}   u^m_{x_k x_j}  )_{x_i}  = S_3 + S_4    \ea     \eeq
 Using  (2.9) with $  k $  playing the role of  $ j $   and  $ x_j $ the role of $ \la $ we have
 \beq  \tag{2.17}   S_4    =  -   p {\ds  \sum_{ i =1}^n   \sum_{ l =1}^N }       u^l_{x_j }
  (  h^{p-2}  u_{ x_j x_i}^l   )_{x_i} = - \De ( h^p)   +  p {\ds  \sum_{ i,j =1}^n   \sum_{ l =1}^N }      h^{p-2}  ( u_{ x_j x_i}^l   )^2  \, .        \eeq
Combining  (2.15)-(2.17) it follows that  $ L (h^p) \geq  p\ \min (1, p-1)  h^{p-2} | \nabla h |^2 $ and thereupon from
   (2.13), (2.14)   that for some $ c = c ( p ) \geq 1, $
  \beq \tag{2.18}   L (h^q )  \geq c^{-1}  q^2   h^{q-2}  | \nabla h |^2  = (4c)^{-1} | \nabla (h^{q/2}) |^2 . \eeq

  Let $ \ph \in C_0^\infty ( B (0, 3/2 ) ) $ and put $ w  = \ph^2. $
 From (2.18), (2.3), (2.7),    and integration  by parts we find for $ q \geq p $ a constant $ c = c ( p, n) \geq 1$ with
 \beq \bea{ll} \tag{2.19}
  T & =  {\ds c ^{-1} \int_{\Om \cap B(0,4)}\phi^2 } |\nabla  h^{q/2}|^2 dx  \leq    {\ds  \,  \int_{\Om \cap B(0,4)} \phi^2  L  (h^{q}) dx} \\ \\ &        {\ds  =   \sum_{i,j=1}^n  \int_{\ar \Om} \phi^2  c_{ij} (h^{q})_{x_j}  \nu_i \, d\si }   -    {\ds  2  \sum_{i,j=1}^n \int_{\Om \cap B(0,4)}  \phi \phi_{x_i} c_{ij} (h^{q})_{x_j} dx  = T_1 + T_2. } \ea  \eeq
From (2.5), smoothness of  $ \ar \Om$  and $u$,  we deduce that $ u_{\nu }  = 0 $ on  $ \ar \Om \cap B ( 0, 2). $ Using this deduction and  (2.12), we see that only the $ \de_{ij}  $ term  contributes  to the  sum  defining $ T_1.  $ This remark and [G]( pp. 133-137 and  Equation 3.1.1.8)   yield

\beq \bea{ll} \tag {2.20}
 T_1= & q {\ds  \int_{\ar \Om} \ph^2 h^{q-2} \sum_{j, k=1}^n \sum_{l=1}^N } u^{l}_{x_k} \nu_j  u^{l}_{x_k x_j}  d \si
 \\ \\ & = {\ds \int_{ \ar \Om } \sum_{l=1}^N  }\ph^2  \, h^{q-2}  \, M ( \nabla_t u^l, \nabla_t u^l )
 \, d \si \leq 0
  \ea \eeq
 where  $M(. , )$ is the second fundamental quadratic  form on $\partial \Om$ and $\nabla_t  u$ is the tangential component of $\nabla u$ on $\partial \Om$. Since $\Om$ is convex, $M(.,)$ is nonpositive. Thus $T_1 \leq 0 $ .
Also using Cauchy's inequality with $ \de$'s  we  get
 \beq \tag{2.21}
 |T_2 |\leq  \delta \int \phi^2 |\nabla (h^{q/2})|^2 \, dx \, +  \, \frac{c }{\delta} \int |\nabla \phi|^2 h^{q} \, dx
\eeq
  Choosing $\delta$ so small that the first term on the righthand side is  $ \leq T/2 $  and then
 using  (2.21),  (2.20),  in  (2.19) we conclude after some arithmetic that \beq \tag{2.22}  \int \phi^2 |\nabla h^{q/2}|^2 dx \leq c   \int |\nabla \phi|^2 h^{q} dx
\end{equation} where $ c = c ( p, n ) \geq 1. $ From (2.22) and Sobolev's inequality for $ \Om \cap B ( 0, 2 ) $ applied to  $ \he  =  h^{q/2}  \ph $  we deduce for $ n \not =  2 $ that
\beq  \tag{2.23}
\begin{split}
\||\he ||^2_{L^{2n/(n-2)} ( \Om \cap B(0,2) ) }\leq (c/|\Om \cap B(0,2) |^2) \,  \| \nabla \he  ||^2_{L^2 (\Om \cap B(0,2))} &\\
\begin{align}
& \leq  ( c /|\Om \cap B(0,2)|^2) \,  \int_{\Om \cap B(0,2) }    |\nabla \ph |^2   h^{q} \, dx
\end{align}
\end{split}
\eeq  If $ n = 2$ replace $ 2n/(n-2) $ by 4 in (2.23).
         We can  now use   Moser iteration (see[GT, ch 8]) in a well known way  and  get for every $ x \in \Om \cap  B (0, 1 ) $  that
\beq \tag{2.24}   ( \ep + |\nabla u |^2 )^{p/2} (x)   \leq  C
  \int_{ \Om \cap B ( 0, 2 ) }  \,  ( \ep +  | \nabla u |^2 )^{p/2} dx \eeq
where  $ C $ has the same dependence as in Theorem 1.1.  Finally,
we  note (see [D], [L], {T] for $N=1$  and [T1] for $N>1$) that   $  \nabla u ( \cdot, \ep )  $ is  H\"{o}lder continuous on compact subsets of  $ \ti \Om \cap B (0, \rho ) $ with H\"{o}lder constants  independent of  $ \ep. $ Also from   (2.6) we find that  $ \{ u ( \cdot, \ep) \} $ is uniformly bounded in $ W^{1,p} (\ti \Om \cap B (0, \rho )). $ Using these facts it  follows  easily  that subsequences of  $ \{ u ( \cdot, \ep) \},  \{ \nabla u ( \cdot, \ep) \} $ converge  uniformly on compact subsets of  $ \Om \cap B ( 0, \rho ) $
as  $ \ep \rar 0  $ to  $ u'  \in  W^{1,p} ( \ti \Om \cap B (0, \rho ) ). $  From uniform integrability type arguments we see that  (2.5)  holds for $ u', \ti \Om \cap B ( 0, \rho ) $ with $\ep =0 $ when  $ v  $ is infinitely differentiable  on  $ \rn{n}. $  Since
a given $ v \in W^{1,p} ( \Om \cap B (0, \rho ) ) $ can be approximated  arbitrarily closely in  the norm of this space by such functions, we conclude that (2.5) is valid  with $ \ep = 0 $ and $ u, \Om $  replaced by
$ u',  \ti \Om. $  Now $ \ti u $ is the unique function (up to a constant) having these properties  so $ u' - \ti u  = $ constant. Thus  (2.6), (2.24)  hold with $ \ep  = 0 $  when $ u, \Om $  are  replaced  by  $ \ti u, \ti \Om$. Hence Theorem 1.1 is true.    $ \Box $
\section{Appendix: Proof of 2.7}
After a rotation if necessary, around a neighborhood of any boundary point $y_0 \in \partial \Om \cap B(0,2)$,  $\partial  \Om$  can be locally represented as $ \{ y', y_n): y_n = \phi(y')\}$  and  $\Om= \{ (y',y_n): y_n > \phi(y') \}$ for some $\phi \in C^{\infty}(R^{n-1})$.  As in  [GT] (Section 14.6), it follows that there exists $ \mu > 0 $ such that if  $ x \in \Om $  and   $ d ( x, \ar \Om )  = d ( x )  <  \mu, $ then  there is a unique point $y(x) \in \partial \Om$, such that $|y-x|=d(x, \partial \Om)= d(x)$. The points $x$ and $y$ are related by

\beq \tag{3.1}
x= y -  \nu(y) d
\eeq
where $\nu$ is the outer unit  normal to $\Om$

\medskip

As proved in Lemma 14.16 in [GT],  the map $ g(y', d)= (y', \phi (y')) - \nu(y', \phi(y')) d$ is locally invertible in a neighborhood $U$ of $(y'_0, 0)$  and maps $U \cap \{ d > 0 \}$ into $\Om.$   Hence by the inverse function theorem
    $y$ and $d$ in (3.1)  are  locally $  C^\infty $ functions of $x$.  Using    $ \nabla d ( x ) = - \nu ( y )  $ one now calculates that  the inverse of  $ g $ is the mapping,
    \beq \tag{3.2}  g^{-1} ( x ) =   ( y', d )  \mbox{ where  $ d = d ( x ) $ and }  y_i =   x_i  - d ( x ) d_{x_i}
    (x)  , 1 \leq i \leq n - 1. \eeq
   By shrinking $U$ if necessary, we may assume that $ g(U \cap \{d>0\})$ is contained in $\Om \cap B(0,5/2)$.

We note  from (3.2) that at $ x \in g ( U ), $
\beq \tag{3.3}
\lan \nabla y_i , \nabla d \ran= \lan e_i - ( \nabla d_{x_i} )  d -  d_{x_i} \nabla d  , \nabla d \ran
= - d \lan  \nabla d_{x_i} , \nabla d \ran = 0
\eeq
where  in the last inequality we have used the fact that  $   \nabla d $ is constant  along the line from
$  x $ to  $ y $ in (3.1) (so the directional derivative of  $ d_{x_i}$ in the direction of  $ \nabla d (x) $  =  0).

 Define $\hat{u}(y',d)= u(x)$ where $u$ is as in (2.5).  We have
\beq \tag{3.4}
\lan \nabla u, \nabla v \ran= \lan \sum_{i=1}^{n-1} \hat{u}_{y_i} \nabla y_i + \hat{u}_{d} \nabla d,  \sum_{i=1}^{n-1} \hat{v}_{y_i} \nabla y_i + \hat{v}_{d} \nabla d \ran
\eeq
Using the orthogonality of $\nabla y_i$ and $\nabla d$, we get
\beq \tag{3.5}
\lan \nabla u, \nabla v \ran =  \hat{u}_{d} \hat{v}_{d}  + \sum_{i,j=1}^{n-1} \hat{u}_{y_i} \hat{v}_{y_j}  q_{ij}(y',d)
\eeq
In conclusion, by renaming the variable $d= y_n$ and using  the change of variables formula,
we deduce from (2.8), (3.5)  that $\hat{u}$  satisfies
\beq \tag{3.6}
{\ds \int_{U \cap y_n > 0 }  (\ep + <A(y)\nabla \hat{u}, \nabla \hat{u}>)^{p/2-1} < A(y) \nabla \hat{u},
\nabla \zeta> V(y',y_n) dy=0}
\eeq
where $U$ is as above and $\zeta$ is in $W^{1,p}_{0} (U)$. $V$ is the Jacobian of  $ g. $
The matrix $A= (a_{ij})$  is defined by   $a_{ij}= q_{ij}$ when $1\leq i,j \leq n-1,$ where $q_{ij}$ is as in (3.5),   $a_{kn}= a_{nk}=0$ when $k<n$ and $a_{nn}=1$.
Also  if $ z \in \rn{n}, $ then $ A ( y ) z  = \eta  $ where
  $ \eta  = (\eta_1, \dots, \eta_n ) $ with
  \[       \eta_i =  \sum_{j=1}^n  a_{ij} ( y ) z_j \,  \mbox{ and }    A  \nabla \hat u = ( A \nabla \hat{u}^1, \dots, A \nabla \hat u^N ).  \] Now let $\xi= \zeta V$. Using $\xi$ in place of $\zeta$ in (3.6), we obtain

\beq \tag{3.7}
\begin{split}
{\ds \int_{U \cap \{ y_n > 0 \}}  (\ep + <A(y)\nabla \hat{u}, \nabla \hat{u}>)^{p/2-1} <A(y)\nabla \hat{u}, \nabla \xi > dy } &  \\
\begin{align}
& = {\ds \int_{U \cap \{ y_n>0\}} (\ep + <A(y) \nabla \hat{u}, \nabla \hat{u}>)^{p/2-1}\sum_{l=1}^N < A(y)
\nabla \hat{u}^l, \nabla  V>\xi^{l} V^{-1} dy}
\end{align}
\end{split}
\eeq
Now we  extend  $A$, $\hat{u}$ and $V$ to $ U \cap \{y_n < 0 \}$  by even reflection and denote them by $B$, $\hat{v}$ and $W$ respectively.  That is,  $B (y', y_n) = A(y',-y_n), \hat v (y', y_n ) = \hat u (y', -y_n), $ and $  W(y',y_n)=V(y',-y_n)$ for $y_n < 0 $. We note that $B$ and $W$ are Lipschitz extensions of $A$ and $V$ respectively and $ \hat v \in  W^{1,p}(U)$

\medskip

Now for any $\zeta   \in W^{1,p}_{0}(U)$
\beq \tag{3.8}
 \int_{U }  (\ep + <B(y) \nabla \hat{v}, \nabla \hat{v}>)^{p/2-1} <B(y)\nabla \hat{v}, \nabla \ze  > dy =  \int_{U \cap \{y_n > 0\} }   + \int_{U \cap \{y_n < 0\} }.
\eeq
For the first integral on the right hand side of (3.8), we use (3.7) and the fact that  $\hat{v},$  $W$ restricted to $\{y_n >0\}$ equal $\hat{u},$  $V$ respectively. For the second integral, we  define $\psi(y', -y_n)= \zeta  (y', y_n)$ and  change variables.   Using the definition of $\hat{v}, \psi$ and $W$  we obtain    (3.7) with $\xi $ replaced by $\psi$.  Altogether we  see that
\beq \tag{3.9}
\begin{split}
{\ds \int_{U}  (\ep + <B(y) \nabla \hat{v},  \nabla \hat{v}>)^{p/2-1} <B(y) \nabla \hat{v}, \nabla \zeta> dy } &  \\
\begin{align}
& = {\ds \int_{U} (\ep + <B(y) \nabla \hat{v},  \nabla \hat{v}>)^{p/2-1}\sum_{l=1}^N < B(y) \nabla \hat{v}^l,
\nabla  W>\zeta^{l} W^{-1} dy}
\end{align}
\end{split}
\eeq
Thus $  \hat v  \in W^{1,p} ( U )  $  is a  weak solution in $ U $ to
 a system of the form
 \beq \tag{3.10}   \nabla \cdot \hat b ( y, \nabla \hat v ) -  \hat b_0 ( y, \nabla \hat v )  = 0.  \eeq
 Moreover this system satisfies  the structural assumption in  [T1] (see (1.1) and (1.7-(1.13) in this paper) except
that in [T1],  the analogue of  $ \hat b $   is assumed to be $ C^1 $ in $ y. $  However  estimates in [T1] only use Lipschitz norms in   $ y $ so are  also valid in our case.
To check the lower order term observe that the $ l $ th component of  $ \hat b_0 ( y, \nabla \hat v ), 1 \leq l \leq N, $ is given by
\beq \tag{3.11}    ( \ep +  <B(y) \nabla \hat v , \nabla  \hat v>)^{p/2-1} < B(y) \nabla  \hat v^l, \nabla  W> W^{-1}.  \eeq
This term may be discontinuous in $  y $  when $ y_n = 0 $ thanks to  $  W_{y_n},   $ but still satisfies the growth and structure assumptions in (1.13) of  [T1].  Therefore, we conclude from [T1], that $\hat{v}  \in  C^{1,\alpha}_{loc} (U) \cap W^{2,2}_{loc} (U)$ (for second derivative estimates, see section 4 in [T1]).

   We note from the  above results and  (3.10), (3.11),   that  $  \hat b_0 ( y, \nabla \hat v (y) ) $ is
   H\"{o}lder continuous in $ U $
      since $ \hat v_{y_n} = 0 $ at points in $ U $ where $ y_n = 0. $. Moreover $(b_{ij})_{y_k}$ are Lipschitz continuous when $k < n$.  Therefore, by  using a  difference quotient argument we
   deduce that if  $ \la = y_{k} $ for $ k < n $   then for $ 1 \leq m \leq N, $
                (3.10) can be differentiated with respect to $ \la $ as in the derivation of (2.9) in order to obtain that $ \hat v^m_\la = \lan \nabla \hat v^m, \la  \ran, 1 \leq m \leq N,  $   is  a   weak  solution to  a  uniformly elliptic linear system in $ N $ equations of the form,
    \beq \tag{3.12}    \sum_{m=1}^N  \sum_{i, j=1}^n   \frac{\ar}{\ar y_i} \left(   C^{m l}_{i j} ( y ) \, (\hat v^m_\la)_{y_j}  \right)  =
     \sum_{i=1}^n   (f_i^l)_{y_i} (y)     \eeq   for fixed $ l, $     where $ f_i^l,  C_{ij}^{m l},   1\leq i, j  \leq n,
        1 \leq l, m \leq N, $  are  H\"{o}lder continuous.    We can now apply Theorem 2.2 in Chapter 3 of  [Gi] to conclude that $ \hat v_\la $  has  H\"{o}lder continuous derivatives  in $ U $ ,  i.e $\hat{v}_{y_i y_j}$  is H\"{o}lder continuous  in $U$ when $ i + j < 2n$. Now  coupled with the fact that  $\hat v$ is in $W^{2,2} ( U )$  we can write (3.10) in non-divergence form and obtain for each $l= 1, ........N$  that

\beq \tag{3.13}
 \sum_{m=1}^N a^{lm} \hat u^{m}_{y_n y_n}  = h^{l}  \mbox{ at $ y \in U $ with $ y_n > 0 $ }
\eeq
where $h^{l}$ is H\"{o}lder continuous in the closure of  $ U  \cap  \{ y  : y_n  > 0 \}   $ and
\[
 a^{lm}=  \delta_{lm} + (p-2)\frac{ \hat u^{l}_{y_n} \hat u^{m}_{y_n}} { ( \ep +  < A(x)\nabla \hat{u}, \nabla \hat{u}>)}
\]  

Thus we see that $ (a^{lm})$ as a matrix is  Holder continuous and positive definite at each point. Consequently, the linear equations corresponding to (3.13) can be solved and $ \hat u^{l}_{y_n y_n}$ is expressible in terms of functions which are H\"{o}lder continuous in the closure of 
$ U  \cap   \{ y_n > 0 \}. $   Therefore $\hat{u} \in C^{2,\alpha}$ upto $y_n=0$. Using the definition of $\hat{u}$,   we then  obtain   that $u$ is $C^{2,\alpha}$  in  the closure of    $ \Om \cap g ( U ).  $
Interior estimates are  similar so we conclude that  $ u  $  is  $ C^{2, \alpha} $  for some $ \al > 0 $ in the closure of  
$ \Om  \cap B ( 0, 2 ). $   
\medskip

         We can now use Schauder type arguments as in
         [ADN] for equations and [ADN1] for systems (see also [Li] for quasilinear equations)  to bootstrap
                   and eventually deduce that $u$ is infinitely differentiable in  the closure of $\Om \cap B(0.2)$. The proof of (2.7) is now complete.\\
 \\

  \noindent \emph{Remark: We emphasize  that in order to satisfy the hypotheses in [T1]  (i.e to get an  even  Lipschitz   extension of  $ A), $ it  was  important that  $a_{kn}=0$ for $k<n$.     This is precisely why  we chose our coordinates  using the  distance function.}

\end{document}